\documentclass[12pt]{amsart} \textwidth=14.5cm \oddsidemargin=1cm
\usepackage{amsmath}
\usepackage{amsxtra}
\usepackage{amscd}
\usepackage{amsthm}
\usepackage{amsfonts}
\usepackage{amssymb}
\usepackage{eucal}

\input prepictex
\input pictex
\input postpictex

\newtheorem{thm}{Theorem}[section]
\newtheorem{lem}{Lemma}[section]

\newtheorem{prop}{Proposition}[section]

\theoremstyle{definition}

\theoremstyle{remark}
\newtheorem{rem}{Remark}[section]

\numberwithin{equation}{section}

\begin{document}

\newcommand{\thmref}[1]{Theorem~\ref{#1}}
\newcommand{\secref}[1]{Section~\ref{#1}}
\newcommand{\lemref}[1]{Lemma~\ref{#1}}
\newcommand{\propref}[1]{Proposition~\ref{#1}}
\newcommand{\corref}[1]{Corollary~\ref{#1}}
\newcommand{\remref}[1]{Remark~\ref{#1}}
\newcommand{\eqnref}[1]{(\ref{#1})}
\newcommand{\exref}[1]{Example~\ref{#1}}

\newcommand{\nc}{\newcommand}
\nc{\on}{\operatorname} \nc{\Z}{{\mathbb Z}} \nc{\C}{{\mathbb C}}
\nc{\oo}{{\mf O}} \nc{\N}{{\mathbb N}} \nc{\bib}{\bibitem}
\nc{\pa}{\partial} \nc{\F}{{\mf F}} \nc{\rarr}{\rightarrow}
\nc{\larr}{\longrightarrow} \nc{\al}{\alpha} \nc{\ri}{\rangle}
\nc{\lef}{\langle} \nc{\W}{{\mc W}} \nc{\gam}{\ol{\gamma}}
\nc{\Q}{\ol{Q}} \nc{\q}{\widetilde{Q}} \nc{\la}{\lambda}
\nc{\ep}{\epsilon} \nc{\g}{\mf g} \nc{\h}{\mf h} \nc{\n}{\mf n}
\nc{\A}{{\mf a}} \nc{\G}{{\mf g}} \nc{\D}{\mc D} \nc{\Li}{{\mc L}}
\nc{\La}{\Lambda} \nc{\is}{{\mathbf i}} \nc{\V}{\mf V}
\nc{\bi}{\bibitem} \nc{\NS}{\mf N}
\nc{\dt}{\mathord{\hbox{${\frac{d}{d t}}$}}} \nc{\E}{\mc E}
\nc{\ba}{\tilde{\pa}} \nc{\half}{\frac{1}{2}}
\def\smapdown#1{\big\downarrow\rlap{$\vcenter{\hbox{$\scriptstyle#1$}}$}}
\nc{\mc}{\mathcal} \nc{\mf}{\mathfrak} \nc{\ol}{\fracline}
\nc{\el}{\ell} \nc{\etabf}{{\bf \eta}} \nc{\x}{{\bf x}}
\nc{\xibf}{{\bf \xi}} \nc{\y}{{\bf y}} \nc{\WW}{\mc W}
\nc{\SW}{\mc S \mc W} \nc{\sd}{\mc S \mc D} \nc{\hsd}{\widehat{\mc
S\mc D}} \nc{\parth}{\partial_{\theta}} \nc{\cwo}{\C[w]^{(1)}}
\nc{\cwe}{\C[w]^{(0)}} \nc{\hf}{\frac{1}{2}}
\nc{\hsdzero}{{}^0\widehat{\sd}} \nc{\hsdpp}{{}^{++}\widehat{\sd}}
\nc{\hsdpm}{{}^{+-}\widehat{\sd}}
\nc{\hsdmp}{{}^{-+}\widehat{\sd}}
\nc{\hsdmm}{{}^{--}\widehat{\sd}} \nc{\gltwo}{{\rm
gl}_{\infty|\infty}} \nc{\btwo}{{B}_{\infty|\infty}}
\nc{\htwo}{{\h}_{\infty|\infty}} \nc{\hglone}{\widehat{\rm
gl}_{\infty}} \nc{\hgltwo}{\widehat{\rm gl}_{\infty|\infty}}
\nc{\hbtwo}{\hat{B}_{\infty|\infty}}
\nc{\hhtwo}{\hat{\h}_{\infty|\infty}} \nc{\glone}{{\rm gl}_\infty}
\nc{\gl}{{\rm gl}} \nc{\ospd}{\mc B} \nc{\hospd}{\widehat{\mc B}}
\nc{\pd}{\mc P} \nc{\hpd}{\widehat{\pd}} \nc{\co}{\mc O}
\nc{\Oe}{\co^{(0)}} \nc{\Oo}{\co^{(1)}} \nc{\sdzero}{{}^0{\sd}}
\nc{\hz}{\hf+\Z} \nc{\vac}{|0 \rangle}

\advance\headheight by 2pt

\title{Infinite-dimensional Lie superalgebras and hook Schur functions}

\author[Shun-Jen Cheng]{Shun-Jen Cheng$^*$}
\thanks{$^*$Partially supported by NSC-grant 91-2115-M-002-007 of the R.O.C}
\address{Department of Mathematics, National Taiwan University, Taipei,
Taiwan 106} \email{chengsj@math.ntu.edu.tw}

\author[Ngau Lam]{Ngau Lam$^{**}$}
\thanks{$^{**}$Partially supported by NSC-grant 90-2115-M-006-015 of the
R.O.C}
\address{Department of Mathematics, National Cheng-Kung
University, Tainan, Taiwan 701} \email{nlam@mail.ncku.edu.tw}

\begin{abstract}
Making use of a Howe duality involving the infinite-dimensional
Lie superalgebra $\hgltwo$ and the finite-dimensional group $GL_l$
of \cite{CW3} we derive a character formula for a certain class of
irreducible quasi-finite representations of $\hgltwo$ in terms of
hook Schur functions. We use the reduction procedure of $\hgltwo$
to $\hat{gl}_{n|n}$ to derive a character formula for a certain
class of level $1$ highest weight irreducible representations of
$\hat{gl}_{n|n}$, the affine Lie superalgebra associated to the
finite-dimensional Lie superalgebra $gl_{n|n}$. These modules turn
out to form the complete set of integrable
$\hat{gl}_{n|n}$-modules of level $1$.  We also show that the
characters of all integrable level $1$ highest weight irreducible
$\hat{gl}_{m|n}$-modules may be written as a sum of products of
hook Schur functions.

\vspace{.3cm}




\end{abstract}

\maketitle

\section{Introduction}

Symmetric functions have been playing an important role in
relating combinatorics and representation theory of Lie
groups/algebras. Interesting combinatorial identities involving
symmetric functions, more often than not, have remarkable
underlying representation-theoretic explanation. As an example
consider the classical Cauchy identity
\begin{equation}\label{cauchy1}
\prod_{i,j}\frac{1}{(1-x_iy_j)} =\sum_{\la} s_\la(x_1,x_2,\cdots)
s_\la(y_1,y_2,\cdots),
\end{equation}
where $x_1,x_2,\cdots$ and $y_1,y_2,\cdots$ are indeterminates,
and $s_\la(x_1,x_2,\cdots)$ stands for the Schur function
associated to the partition $\la$.  Here the summation of $\la$
above is over all partitions. Now the underlying
representation-theoretic interpretation of \eqnref{cauchy1} is of
course the so-called $(GL,GL)$ Howe duality \cite{H1} \cite{H2}.
Namely, let $\C^m$ and $\C^n$ be the $m$- and $n$-dimensional
complex vector spaces, respectively.  We have an action of the
respective general linear groups $GL_m$ and $GL_n$ on $\C^m$ and
$\C^n$. This induces a joint action of $GL_m\times GL_n$ on
$\C^m\otimes\C^n$, which in turn induces an action on the
symmetric tensor $S(\C^m\otimes\C^n)$. As partitions of
appropriate length may be regarded as highest weights of
irreducible representations of a general linear group,
\eqnref{cauchy1} simply gives an identity of characters of the
decomposition of $S(\C^m\otimes\C^n)$ with respect to this joint
action.

It was observed in \cite{BR} that a generalization of Schur
functions, the so-called the hook Schur functions (see
\eqnref{hookschur1} for definition), play a similar role in the
representation theory of a certain class of finite-dimensional
irreducible modules over the general linear Lie superalgebra. To
be more precise, consider the general linear Lie superalgebra
$gl_{m|n}$ acting on the complex superspace $\C^{m|n}$ of
(super)dimension $(m|n)$. We may consider its induced action on
the $k$-th tensor power $T^k(\C^{m|n})=\bigotimes^k(\C^{m|n})$. It
turns out \cite{BR} that the tensor algebra
$T(\C^{m|n})=\bigoplus_{k=0}^\infty T^k(\C^{m|n})$ is completely
reducible as a $gl_{m|n}$-module and the characters of the
irreducible representations appearing in this decomposition are
given by hook Schur functions associated to partitions lying in a
certain hook whose shape is determined by the integers $m$ and
$n$. Now, as in the classical case, one may consider the joint
action of two general linear Lie superalgebras $gl_{m|n}\times
gl_{p|q}$ on the symmetric tensor $S(\C^{m|n}\otimes\C^{p|q})$.
This action again is completely reducible and its decomposition
with respect to the joint action, in a similar fashion, gives rise
to a combinatorial identity involving hook Schur functions
\cite{CW1}. So here we have an interplay between combinatorics and
representation theory of finite-dimensional Lie superalgebras as
well.  For another interplay involving Schur $Q$-functions and the
queer Lie superalgebra see \cite{CW2}.  For further articles
related to Howe duality in the Lie superalgebra settings we refer
to \cite {S1}, \cite{S2}, \cite{N} and \cite{OP}.

The purpose of the present paper is to demonstrate that symmetric
functions may play a similarly prominent role relating
combinatorics and representation theory of infinite-dimensional
Lie superalgebras as well. It was shown in \cite{CW3} that on the
infinite-dimensional Fock space generated by $n$ pairs of free
bosons and $n$ pairs of free fermions we have a natural commuting
action of the finite-dimensional group $GL_n$ and the
infinite-dimensional Lie superalgebra $\hgltwo$ of central charge
$n$.  It can be shown \cite{CW3} that the pair $(GL_n,\hgltwo)$
forms a dual pair in the sense of Howe. The irreducible
representations of $GL_n$ appearing in this decomposition ranges
over all rational representations, so that they are parameterized
by generalized partitions of length not exceeding $n$.  The
irreducible representations of $\hgltwo$ appearing in the same
decomposition are certain quasi-finite highest weight irreducible
representations. In particular this relates rational
representations of the finite-dimensional group $GL_n$ and a
certain class of quasi-finite representations of the
infinite-dimensional Lie superalgebra $\hgltwo$.

A natural question that arises is the computation of the character
of these quasi-finite highest weight irreducible representations
of $\hgltwo$. This is solved in the present paper by combining the
Howe duality of \cite{CW3} and a combinatorial identity involving
hook Schur functions (\propref{combid}). It turns out that the
characters of these representations can be written as an infinite
sum of products of two hook Schur functions. Each coefficient of
these products can be determined by decomposing a certain tensor
product of two finite-dimensional irreducible representations of
$GL_n$. The same method applied to the classical Lie algebra
$\hat{gl}_\infty$ (using \lemref{independence} now together with
the dual Cauchy identity instead of \propref{combid}) gives rise
to a character formula for $\hat{gl}_\infty$ involving Schur
functions instead of hook Schur functions. This formula has been
discovered earlier by Kac and Radul \cite{KR2} and even earlier in
the simplest case by Awata, Fukuma, Matsuo and Odake \cite{AFMO2}.
However, our approach in these classical cases appears to be
simpler. It is indeed remarkable that the same character identity
obtained in \cite{KR2} with Schur functions replaced by hook Schur
functions associated to the same partitions gives rise to our
character identity for $\hgltwo$. That is, the coefficients remain
unchanged!

Making use of the same combinatorial identity we then proceed to
compute the corresponding $q$-character formula for this class of
highest weight irreducible representations of $\hgltwo$. We remark
that when computing just the $q$-characters we can obtain a
simpler formula, which involves a sum of just hook Schur
functions, instead of a product of hook Schur functions.  We note
that a $q$-character formula in the case when the central charge
is $1$ has been obtained earlier by Kac and van de Leur \cite{KL}.
Our formula in this case looks rather different from theirs, thus
giving rise to another combinatorial identity.  It would be
interesting to find a purely combinatorial proof of this identity.

We use the reduction procedure from $\hgltwo$ to the
$\hat{gl}_{n|n}$ \cite{KL} in the level $1$ case to obtain a
character formula for certain highest weight irreducible
representations of the affine Lie superalgebra $\hat{gl}_{n|n}$ at
level $1$.  However, the Borel subalgebra coming from the
reduction procedure is different from the standard Borel
subalgebra, and hence the corresponding highest weights in general
are different.  However, we show that by a sequence of odd
reflections \cite{PS} our highest weights may be transformed into
highest weights corresponding to integrable
$\hat{gl}_{n|n}$-modules (in the sense of \cite{KW2}) so that we
obtain a character formula for all integrable highest weight
$\hat{gl}_{n|n}$-modules of level $1$. In \cite{KW2} a character
formula has been obtained for level $1$ integrable highest weight
irreducible $\hat{gl}_{m|n}$-modules. Our formula looks rather
different.

We also show that by applying our method together with \cite{KW2}
the characters of all level $1$ integrable representations of
$\hat{gl}_{m|n}$, $m\ge 2$, may be written in terms of hook Schur
functions as well. This seems to indicate the relevance of these
generalized symmetric functions in the representation theory of
affine superalgebras.

As we have obtained a character formula for certain
representations of $\hgltwo$ at arbitrary positive integral level,
it is our hope that our formula may provide some direction in
finding a character formula for integrable $\hat{gl}_{n|n}$-, or
maybe even $\hat{gl}_{m|n}$-modules, at higher positive integral
levels.

The paper is organized as follows. In Section 2 we collect the
definitions and notation to be used throughout. In Section 3 we
first prove the combinatorial identity mentioned above and then
use it to write the characters of certain $\hgltwo$-modules in
terms of hook Schur functions.  In Section 4 we calculate a
$q$-character formula for these modules, while in Section 5 we
calculate the characters of the associated affine
$\hat{gl}_{n|n}$-modules. In Section 6 we compute the tensor
product decomposition of two $\hgltwo$-modules. It turns out that
even though such a decomposition involves an infinite number of
irreducible components, each irreducible component appears with a
finite multiplicity.  This multiplicity can be computed via the
usual Littlewood-Richardson rule.

\section{Preliminaries}

Let $\C^{m|n}=\C^{m|0}\oplus\C^{0|n}$ denote the $m|n$-dimensional
superspace. Let $gl_{m|n}$ be the Lie superalgebra of general
linear transformations on the superspace $\C^{m|n}$. Choosing a
basis $\{e_1,\cdots,e_m\}$ for the even subspace $\C^{m|0}$ and a
basis $\{f_1,\cdots,f_n\}$ for the odd subspace $\C^{0|n}$, we may
regard $gl_{m|n}$ as $(m+n)\times(m+n)$ matrices of the form
\begin{equation}\label{supermatrix}
\begin{pmatrix}
E&a\\
b&e
\end{pmatrix},
\end{equation}
where the complex matrices $E$, $a$, $b$ and $e$ are respectively
$m\times m$, $m\times n$, $n\times m$ and $n\times n$.  Let
$X_{ij}$ denote the corresponding elementary matrix with $1$ in
the $i$-th row and $j$-th column and zero elsewehre, where
$X=E,b,a,e$. Then $\h=\sum_{i=1}^{m}\C E_{ii}+\sum_{j=1}^n\C
e_{jj}$ is a Cartan subalgebra of $gl_{m|n}$.

It is clear that any ordering of the basis
$\{e_1,\cdots,e_m,f_1,\cdots,f_n\}$ that preserves the order among
the even and odd basis elements themselves gives rise to a Borel
subalgebra of $gl_{m|n}$ containing $\h$. In particular the
ordering $e_1<\cdots<e_m<f_1<\cdots<f_n$ gives rise to the
standard Borel subalgebra.  In the case when $m=n$ the ordering
$f_1<e_1<f_2<e_2<\cdots<f_n<e_n$ gives rise to a Borel subalgebra
that we will refer to as {\em non-standard} from now on.

Fixing the standard Borel subalgebra, we let $V^\la_{m|n}$ denote
the finite-dimensional highest weight irreducible module with
highest weight $\la$.

Let $\epsilon_i\in\h^*$ be defined by
$\epsilon_i(E_{jj})=\delta_{ij}$ and $\epsilon_i(e_{jj})=0$.
Furthermore let $\delta_j$ be defined by $\delta_j(E_{ii})=0$ and
$\delta_j(e_{ii})=\delta_{ij}$. Then $\epsilon_i$ and $\delta_j$
are the fundamental weights of $gl_{m|n}$.

Let $\C[t,t^{-1}]$ be the ring of Laurent polynomials in the
indeterminate $t$. Let $\hat{gl}_{m|n}\equiv
gl_{m|n}\otimes\C[t,t^{-1}]+\C C+\C d$ be the affine Lie
superalgebra associated to the Lie superalgebra $gl_{m|n}$.
Writing $A(k)$ for $A\otimes t^k$, $A\in gl_{m|n}$, the Lie
(super)bracket is given by
\begin{align*}
&[A(k),B(l)]=[A,B](k+l)+\delta_{k+l,0}k{\rm Str}(AB)C,\\
&[d,A(k)]=k A(k),\quad A,B\in gl_{m|n},\ k,l\in\Z,
\end{align*}
Here $C$ is a central element, $d$ is the scaling element and $\rm
Str$ denote the super trace operator of a matrix, which for a
matrix of the form \eqnref{supermatrix} takes the form ${\rm
Tr}(E)-{\rm Tr}(e)$.

A Cartan subalgebra of $\hat{gl}_{m|n}$ is given by
$\hat{\h}=\h+\C C+\C d$. We may extend respectively $\epsilon_i$
and $\delta_j$ to elements $\tilde{\epsilon}_i$ and
$\tilde{\delta}_j$ in $\hat{\h}^*$ in a trivial way. Furthermore
we define $\tilde{\Lambda}_0\in\hat{\h}^*$ and
$\tilde{\delta}\in\hat{\h}^*$ by
$\tilde{\Lambda}_0(\h)=\tilde{\Lambda}_0(d)=0$,
$\tilde{\Lambda}_0(C)=1$ and
$\tilde{\delta}(\h)=\tilde{\delta}(C)=0$, $\tilde{\delta}(d)=1$,
respectively.

Let $B\subseteq gl_{m|n}$ be a Borel subalgebra containing $\h$.
Then $B+\C C+\C d+ gl_{m|n}\otimes t\C[t]$ is a Borel subalgebra
of $\hat{gl}_{m|n}$. We define highest weight irreducible modules
of $\hat{gl}_{m|n}$ in the usual way.  It is clear that any
highest weight irreducible $\hat{gl}_{m|n}$-module is completely
determined by an element $\Lambda\in\hat{\h}^*$.  We will denote
this module by $L(\hat{gl}_{m|n},\La)$.

Consider now the infinite-dimensional complex superspace
$\C^{\infty|\infty}$ with even basis elements labelled by integers
and odd basis elements labelled half-integers.  Arranging the
basis elements in strictly increasing order any linear
transformation may be written as an infinite-sized square matrix
with coefficients in $\C$.  This associative algebra is naturally
$\Z_2$-graded, so that it is an associative superalgebra, which we
denote by $\tilde{M}_{\infty|\infty}$. Let
\begin{equation*}
M_{\infty|\infty}:=\{A=(a_{ij})\in\tilde{M}_{\infty|\infty},
i,j\in\hf\Z|\ a_{ij}=0\ {\rm for}\ |j-i|>>0\}.
\end{equation*}
That is, $M_{\infty|\infty}$ consists of those matrices in
$\tilde{M}_{\infty|\infty}$ with finitely many non-zero diagonals.
We denote the corresponding Lie superalgebra by
$gl_{\infty|\infty}$.  Furthermore let us denote by $e_{ij}$,
$i,j\in\half\Z$ the elementary matrices with $1$ at the $i$-th row
and  $j$-th column and $0$ elsewhere.  Then the subalgebra
generated by $\{e_{ij}|i,j\in\half\Z\}$ is a dense subalgebra
inside $gl_{\infty|\infty}$.

The Lie superalgebra $gl_{\infty|\infty}$ has a central extension
(by an even central element $C$), denoted from now on by
$\hgltwo$, corresponding to the following two-cocycle
\begin{equation*}
\alpha(A,B)={\rm Str}([J,A]B),\quad A,B\in\gltwo,
\end{equation*}
where $J$ denotes the matrix $\sum_{r\le 0}e_{rr}$, and for a
matrix $D=(d_{ij})\in\gltwo$, $\rm Str(D)$ stands for the
supertrace of the matrix $D$ and which here is given by
$\sum_{r\in\hf\Z}(-1)^{2r}d_{rr}$. We note that the expression
$\alpha(A,B)$ is well-defined for $A,B\in\gltwo$.

The Lie superalgebra $\hat{gl}_{\infty|\infty}$ has a natural
$\half\Z$-gradation by setting ${\rm deg}E_{ij}=j-i$, for
$i,j\in\half\Z$.  Thus we have the triangular decomposition
\begin{equation*}
\hat{gl}_{\infty|\infty}=(\hat{gl}_{\infty|\infty})_{-}\oplus
(\hat{gl}_{\infty|\infty})_0\oplus (\hat{gl}_{\infty|\infty})_+,
\end{equation*}
where the subscripts $+$, $0$ and $-$ respectively denote the
positive, zero-th and negative graded components.  Thus we have a
notion of a highest weight Verma module, which contains a unique
irreducible quotient, which is determined by an element
$\La\in(\hgltwo)_0^*$.  We will denote this module by
$L(\hgltwo,\La)$.  Let $\omega_s$, $s\in\hf\Z$, denote the
fundamental weights of $\hgltwo$. That is, $\omega_s(e_{rr})=0$,
$r\in\hf\Z$, and $\omega_s(C)=0$. Furthermore let
$\Lambda_0\in(\hgltwo)_0^*$ with $\Lambda_0(e_{rr})=0$ and
$\Lambda_0(C)=1$.

Note that by declaring the highest weight vectors to be of degree
zero, the module $L(\hgltwo,\La)$ is naturally $\half\Z$-graded,
i.e.
\begin{equation*}
L(\hgltwo,\La)=\oplus_{r\in\half\Z_+}L(\hgltwo,\La)_r.
\end{equation*}
The module $L(\hgltwo,\La)$ is said to be {\em quasi-finite}
\cite{KR1} if ${\rm dim}L(\hgltwo,\La)_r<\infty$, for all
$r\in\half\Z_+$.

\section{A Character formula for $\hat{gl}_{\infty|\infty}$-modules}

First we recall the notion of the hook Schur function of
Berele-Regev \cite{BR}.

Let $\x=\{x_1,x_2,\cdots\}$ be a countable set of variables. To a
partition $\la$ of non-negative integers we may associate the
Schur function $s_\la(x_1,x_2,\cdots)$.  We will write $s_\la(\x)$
for $s_\la(x_1,x_2,\cdots)$.  For a partition $\mu\subset\la$ we
let $s_{\la/\mu}(\x)$ denote the corresponding skew Schur
function. Denoting by $\mu'$ the conjugate partition of a
partition $\mu$ the {\em hook Schur function} corresponding to a
partition $\la$ is defined by
\begin{equation}\label{hookschur1}
HS_{\la}(\x;\y):=\sum_{\mu\subset\la}s_\mu(\x)s_{\la'/\mu'}(\y),
\end{equation}
where as usual $\y=\{y_1,y_2,\cdots\}$.

Let $\la$ be a partition and $\mu\subseteq\la$.  We fill the boxes
in $\mu$ with entries from the linearly ordered set
$\{x_1<x_2<\cdots\}$ so that the resulting tableau is
semi-standard. Recall that this means that the rows are
non-decreasing, while the columns are strictly increasing.  Next
we fill the skew partition $\la/\mu$ with entries from the
linearly ordered set $\{y_1<y_2<\cdots\}$ so that it is conjugate
semi-standard, which means that the rows are strictly increasing,
while its columns are non-decreasing.  We will refer to such a
tableau as an {\em $(\infty|\infty)$-semi-standard tableau}
(cf.~\cite{BR}).  To each such tableau $T$ we may associate a
polynomial $(\x\y)^T$, which is obtained by taking the products of
all the entries in $T$. Then we have \cite{BR}
\begin{equation}\label{hookschur}
HS_{\la}(\x;\y)=\sum_{T}(\x\y)^T,
\end{equation}
where the summation is over all $(\infty|\infty)$-semi-standard
tableaux of shape $\la$.

We have the following combinatorial identity involving hook Schur
functions that is crucial in the sequel.

\begin{prop}\label{combid}
Let $\x=\{x_1,x_2,\cdots\}$, $\y=\{y_1,y_2,\cdots\}$ be two
infinite countable sets of variables and ${\bf
z}=\{z_1,z_2,\cdots,z_m\}$ be $m$ variables.  Then
\begin{equation}
\prod_{i,j,k}(1-x_iz_k)^{-1}(1+y_jz_k)=\sum_{\la}HS_{\la}(\x;\y)s_{\la}(\bf
z),
\end{equation}
where $1\le i,j<\infty$, $1\le k\le m$ and $\la$ is summed over
all partitions $\la$ with length not exceeding $m$.
\end{prop}

\begin{proof}
Consider the classical Cauchy identity
\begin{equation}\label{cauchy}
\prod_{i,j}(1-x_iz_k)^{-1}(1-y_jz_k)^{-1}=\sum_{\la}s_\la(\x,\y)s_\la({\bf
z}),
\end{equation}
where $\la$ is summed over all partitions of length not exceeding
$m$.  Recall that for any partition $\la$ one has (cf.~\cite{M}
(I.5.9))
\begin{equation}\label{symid}
s_{\la}(\x,\y)=\sum_{\mu\subset\la}s_\mu(\x)s_{\la/\mu}(\y).
\end{equation}
Let $\omega$ denote the involution of the ring of symmetric
functions, which sends the elementary symmetric functions to the
complete symmetric functions, so that we have
$\omega(s_{\la}(\x))=s_{\la'}(\x)$.  Now applying $\omega$ to the
set of variables $\y$ in \eqnref{cauchy} we obtain together with
\eqnref{symid}
\begin{align*}
\prod_{i,j,k}(1-x_iz_k)^{-1}(1+y_jz_k)&=\sum_{\la}(\sum_{\mu\subset\la}
s_\mu(\x)s_{\la'/\mu'}(\y))s_\la({\bf z})\\
&=\sum_{\la}HS_{\la}(\x;\y)s_{\la}(\bf z),
\end{align*}
as required.
\end{proof}

We note that \propref{combid} in the case when the sets of
variables are all finite sets follows from the Howe duality
(\cite{H1}, \cite{H2}) involving a general linear Lie superalgebra
and a general linear Lie algebra described in \cite{CW1}. Since we
will need this result in the case when both algebras involved are
Lie algebras later on we will recall it here.

\begin{prop}\cite{H2}\label{aux31}
The Lie algebras $gl_d$ and $gl_{m}$ with their natural actions on
$S(\C^d\otimes\C^m)$ form a dual pair. With respect to their joint
action we have the following decomposition.
\begin{equation*}
S(\C^d\otimes\C^m)\cong\sum_\la V^\la_d\otimes V^\la_m,
\end{equation*}
where the summation is over all partitions with length not
exceeding ${\rm min}(l,m)$.
\end{prop}

Below we will recall the $\hat{gl}_{\infty|\infty}\times gl_l$
duality of \cite{CW3}.  Consider $l$ pairs of free fermions
$\psi^{\pm,i}(z)$ and $l$ pairs of free bosons $\gamma^{\pm,i}(z)$
with $i=1,\cdots,l$.  That is we have
\begin{align*}
\psi^{+,i}(z)&=\sum_{n\in\Z}\psi^{+,i}_nz^{-n-1},\quad\quad\
\psi^{-,i}(z)=\sum_{n\in\Z}\psi^{-,i}_nz^{-n},\\
\gamma^{+,i}(z)&=\sum_{r\in\frac{1}{2}+\Z}\gamma^{+,i}_rz^{-r-1/2},\quad
\gamma^{-,i}(z)=\sum_{r\in\frac{1}{2}+\Z}\gamma^{-,i}_rz^{-r-1/2}
\end{align*}
with non-trivial commutation relations
$[\psi^{+,i}_m,\psi^{-,j}_n]=\delta_{ij}\delta_{m+n,0}$ and
$[\gamma^{+,i}_r,\gamma^{-,j}_s]=\delta_{ij}\delta_{r+s,0}$.

Let $\F$ denote the corresponding Fock space generated by the
vaccum vector $|0>$.  That is
$\psi^{+,i}_n|0>=\psi^{-,i}_m|0>=\gamma^{\pm,i}_r|0>=0$, for $n\ge
0$, $m>0$ and $r>0$. These operators are called {\em annihilation
operators}.

Explicitly we have an action of $\hat{gl}_{\infty|\infty}$ of
central charge $l$ on $\F$ given by ($i,j\in\Z$ and
$r,s\in\frac{1}{2}+\Z$)
\begin{align}
e_{ij}&=\sum_{p=1}^l:\psi^{+,p}_{-i}\psi^{-,p}_{j}:,\nonumber\\
e_{rs}&=-\sum_{p=1}^l:\gamma^{+,p}_{-r}\gamma^{-,p}_{s}:,\nonumber\\
e_{is}&=\sum_{p=1}^l:\psi^{+,p}_{-i}\gamma^{-,p}_{s}:,\nonumber\\
e_{rj}&=-\sum_{p=1}^l:\gamma^{+,p}_{-r}\psi^{-,p}_{j}:.\nonumber
\end{align}

An action of $gl_l$ on $\F$ is given by the formula
\begin{equation*}
E_{ij}=\sum_{n\in\Z}:\psi^{+,i}_{-n}\psi^{-,j}_{n}:-
\sum_{r\in1/2+\Z}:\gamma^{+,i}_{-r}\gamma^{-,j}_{r}:.
\end{equation*}

Here and further $::$ denotes the normal ordering of operators.
That is, if $A$ and $B$ are two operators, then $:AB:=AB$, if $B$
is an annihilation operator, while $:AB:=(-1)^{{\rm p}(A){\rm
p}(B)}BA$, otherwise.  As usual, ${\rm p}(X)$ denotes the parity
of the operator $X$.

Before stating the duality of \cite{CW3} we need some more
notation. For $j\in\Z_+$ we define the matrices $X^{-j}$ as
follows:
\begin{eqnarray*}
 X^0=&
\begin{pmatrix}
\psi_0^{-,l}&\psi_0^{-,l-1}&\cdots &\psi_0^{-,1}\\
\psi_0^{-,l}&\psi_0^{-,l-1}&\cdots &\psi_0^{-,1}\\
\vdots&\vdots&\cdots &\vdots\\
\psi_0^{-,l}&\psi_0^{-,l-1}&\cdots
&\psi_0^{-,1}\\
\end{pmatrix},\allowdisplaybreaks\\
 X^{-1}=&
\begin{pmatrix}
\gamma_{-\hf}^{-,l}&\gamma_{-\hf}^{-,l-1}&\cdots
&\gamma_{-\hf}^{-,1}\\ \psi_{-1}^{-,l}&\psi_{-1}^{-,l-1}&\cdots
&\psi_{-1}^{-,1}\\ \vdots&\vdots&\cdots &\vdots\\
\psi_{-1}^{-,l}&\psi_{-1}^{-,l-1}&\cdots &\psi_{-1}^{-,1}\\
\end{pmatrix},\allowdisplaybreaks\\
 X^{-2}=&
\begin{pmatrix}
\gamma_{-\hf}^{-,l}&\gamma_{-\hf}^{-,l-1}&\cdots
&\gamma_{-\hf}^{-,1}\\
\gamma_{-\frac{3}{2}}^{-,l}&\gamma_{-\frac{3}{2}}^{-,l-1}&\cdots
&\gamma_{-\frac{3}{2}}^{-,1}\\
\psi_{-2}^{-,l}&\psi_{-2}^{-,l-1}&\cdots &\psi_{-2}^{-,1}\\
\vdots&\vdots&\cdots &\vdots\\
\psi_{-2}^{-,l}&\psi_{-2}^{-,l-1}&\cdots &\psi_{-2}^{-,1}\\
\end{pmatrix},\allowdisplaybreaks\\
&\vdots\\&\vdots\\
 X^{-k} \equiv X^{-l}=&
\begin{pmatrix}
\gamma_{-\hf}^{-,l}&\gamma_{-\hf}^{-,l-1}&\cdots
&\gamma_{-\hf}^{-,1}\\
\gamma_{-\frac{3}{2}}^{-,l}&\gamma_{-\frac{3}{2}}^{-,l-1}&\cdots
&\gamma_{-\frac{3}{2}}^{-,1}\\ \vdots&\vdots&\cdots &\vdots\\
\gamma_{-l+\hf}^{-,l}&\gamma_{-l+\hf}^{-,l-1}&\cdots
&\gamma_{-l+\hf}^{-,1}\\
\gamma_{-l-\hf}^{-,l}&\gamma_{-l-\hf}^{-,l-1}&\cdots
&\gamma_{-l-\hf}^{-,1}\\
\end{pmatrix}, \quad k \ge l.\allowdisplaybreaks\\
\end{eqnarray*}
The matrices $ X^j$, for $j\in\N$, are defined similarly. Namely,
$ X^j$ is obtained from $ X^{-j}$ by replacing $\psi^{-,k}_i$ by
$\psi^{+,l-k+1}_i$ and $\gamma^{-,k}_r$ by $\gamma^{+,l-k+1}_r$.

For $0\le r\le l$, we let $X^i_{r}$($i \ge 0$) denote the first
$r\times r$ minor of the matrix $ X^i$ and let $X^i_{-r}$($i <0$)
denote the first $r\times r$ minor of the matrix $ X^i$.

Consider a generalized partition $\la=(\la_1,\la_2,\cdots,\la_p)$
of length not exceeding $l$ with
\begin{equation*}
\la_1\ge\la_2\ge\cdots\ge\la_i>\la_{i+1}
=0=\cdots=\la_{j-1}>\la_j\ge\cdots\ge\la_l.
\end{equation*}

Now the irreducible rational representations of ${\rm GL}_l$ are
parameterized by generalized partitions, hence these may be
interpreted as highest weights of irreducible representations of
${\rm GL}_l$.  We denote the corresponding finite-dimensional
highest weight irreducible ${\rm GL}_l$- (or $gl_l$-) module by
$V^\la_l$. Let $\la_{j}'$ be the length of the $j$-th column of
$\la$. We use the convention that the first column of $\la$ is the
first column of the partition $\la_1\ge\la_2\ge\cdots\ge\la_i$.
The column to the right is the second column of $\la$, while the
column to the left of it is the zeroth column and the column to
the left of the zeroth column is the $-1$-st column.  We also use
the convention that a non-positive column has non-positive length.
As an example consider $\la=(5,3,2,1,-1,-2)$ with $l(\la)=6$. We
have $\la_{-1}'=-1$, $\la_0'=-2$, $\la_1'=4$ etc.~(see
\eqnref{standard}).
\begin{equation}\label{standard}
{\beginpicture \setcoordinatesystem units <1.5pc,1.5pc> point at 0
2 \setplotarea x from -1.5 to 1.5, y from -2 to 4 \plot 0 0 0 4 1
4 1 0 0 0 / \plot 1 1 2 1 2 4 1 4 / \plot 2 2 3 2 3 4 2 4 / \plot
3 3 5 3 5 4 3 4 / \plot 0 1 1 1 / \plot 0 2 2 2 / \plot 0 3 3 3 /
\plot 4 3 4 4 / \plot 0 0 -1 0 / \plot 0 0 0 -2 / \plot -1 0 -1 -2
/ \plot 0 -2 -2 -2 / \plot 0 -1 -2 -1 / \plot -2 -2 -2 -1 /
\endpicture}
\end{equation}

For $\La\in (\hgltwo)_0^*$, we set $\La_s=\La(e_{ss})$, for
$s\in\hf\Z$. Given a generalized partition $\la$ with $l(\la)\le
l$, we define $\Lambda(\la)\in(\hgltwo)_0^*$ by:

\begin{align*}
 &\La(\la)_i= \langle\la_i'-i \rangle,\quad i\in\N,\\
 &\La(\la)_j= -\langle-\la'_j+j \rangle, \quad j\in -\Z_{+},\\
 &\La(\la)_r=\langle\la_{r+1/2}-(r-1/2) \rangle,\quad r\in\hf+\Z_+,\\
 &\La(\la)_{s}=-\langle -\la_{p+(s+1/2)}+({s-1/2}) \rangle,\quad
 s\in-\hf-\Z_+,\\
&\La(\la)(C)=l.
\end{align*}
Here for an integer $k$ the expression $<k>\equiv k$, if $k> 0$,
and $<k>\equiv 0$, otherwise. We have the following theorem.

\begin{thm} \cite{CW3}\label{duality}
The Lie superalgebra $\hat{gl}_{\infty|\infty}$ and $gl_l$ form a
dual pair on $\F$ in the sense of Howe.  Furthermore we have the
following (multiplicity-free) decomposition of $F$ with respect to
their joint action
\begin{equation*}
\F\cong\sum_{\la}L(\hat{gl}_{\infty|\infty},\Lambda(\la))\otimes
V_l^{\la},
\end{equation*}
where the summation is over all generalized partitions of length
not exceeding $l$. Furthermore, the joint highest weight vector of
the $\la$-component is given by
\begin{equation*}
{\rm det} X^{\la_{l}+1}_{{\la}_{\la_{l}+1}'}\cdots{\rm det}
X^{-1}_{{\la}_{-1}'} \cdot {\rm det} X^0_{\la_0'} \cdot {\rm det}
X^1_{\la_1'} \cdot {\rm det} X^2_{\la_2'}\cdots{\rm det}
X^{\la_1}_{{\la}_{\la_1}'}\vac.
\end{equation*}
\end{thm}

We compute for $i\in\Z$, $r\in\frac{1}{2}+\Z$
\begin{align*}
&[e_{ii},\psi^{+,p}_{-n}]=\delta_{in}\psi^{+,p}_{-n},\\
&[e_{ii},\psi^{-,p}_{-n}]=-\delta_{-in}\psi^{-,p}_{-n},\\
&[e_{rr},\psi^{\pm,p}_{-n}]=[e_{ii},\gamma^{\pm,p}_{-r}]=0,\\
&[e_{rr},\gamma^{+,p}_{-s}]=\delta_{rs}\gamma^{+,p}_{-s},\\
&[e_{rr},\gamma^{-,p}_{-s}]=-\delta_{-rs}\gamma^{-,p}_{-s}.
\end{align*}

Furthermore for $i=1,\cdots,l$ we have
\begin{align*}
&[E_{ii},\psi^{+,p}_{-n}]=\delta_{ip}\psi^{+,p}_{-n},\\
&[E_{ii},\psi^{-,p}_{-n}]=-\delta_{ip}\psi^{-,p}_{-n},\\
&[E_{ii},\gamma^{+,p}_{-r}]=\delta_{ip}\gamma^{+,p}_{-r},\\
&[E_{ii},\gamma^{-,p}_{-r}]=-\delta_{ip}\gamma^{-,p}_{-r}.
\end{align*}

Let $e$ be a formal indeterminate and set for $j\in\Z$, $r\in\hf
+\Z$, $i=1,\cdots,l$
\begin{equation*}
x_i=e^{\epsilon_i},\quad y_j=e^{\omega_j},\quad z_r=e^{\omega_r},
\end{equation*}
where $\epsilon_1,\cdots,\epsilon_l$ and $\omega_s$ are the
respective fundamental weights of $gl_l$ and $\hgltwo$ introduced
earlier. It is easy to see that the character of $\F$, with
respect to the abelian algebra $\sum_{s\in\half\Z}\C
e_{ss}\oplus\sum_{i=1}^l\C E_{ii}$, is given by
\begin{equation}\label{character1}
{\rm ch}\F=\prod_{i=1}^l\frac{\prod_{n\in\N}(1+x_iy_n)
\prod_{m\in\Z_+}(1+x_i^{-1}y_{-m}^{-1})}{\prod_{r\in1/2+\Z_+}
(1-x_iz_r)(1-x_i^{-1}z_{-r}^{-1})}
\end{equation}
By \propref{combid} we can rewrite \eqnref{character1} as
\begin{equation}\label{character2}
{\rm ch}\F=\sum_{\mu,\nu}HS_\mu({\bf z};\y)HS_\nu({\bf
z}^{-1};{\bf y^{-1}})s_\mu(\x)s_\nu(\x^{-1}),
\end{equation}
where $\mu$ and $\nu$ are summed over all partitions of length not
exceeding $l$.  Here we use the notation $\y=\{y_1,y_2,\cdots\}$,
$\y^{-1}=\{y_0^{-1},y_{-1}^{-1},\cdots\}$, ${\bf
z}=\{z_{\frac{1}{2}},z_{\frac{3}{2}},\cdots\}$, ${\bf
z}^{-1}=\{z^{-1}_{-\frac{1}{2}},z^{-1}_{-\frac{3}{2}},\cdots\}$,
$\x=\{x_1,x_2,\cdots,x_l\}$, and
$\x^{-1}=\{x^{-1}_1,x_2^{-1},\cdots,x_l^{-1}\}$.

It is clear that $s_{\nu}(\x^{-1})$ is just the character of the
$gl_l$-module $(V_{l}^\nu)^*$, the module contragredient to
$V^\nu_l$. Therefore we have
\begin{equation}\label{branching}
s_\mu(\x)s_\nu(\x^{-1})=\sum_{\la}c^\la_{\mu\nu}{\rm ch}V_{l}^\la,
\end{equation}
where the summation now is over generalized partitions of length
not exceeding $l$.  Here the non-negative integers
$c^\la_{\mu\nu}$ are of course just the multiplicity of $V^\la_l$
in the tensor product decomposition of $V^\mu_l\otimes (V^\nu_l)^*
$. This combined with \eqnref{character2} allows us to write the
character of $\F$ as
\begin{equation}\label{character3}
{\rm ch}\F= \sum_{\la}\Bigl{(}\sum_{\mu,\nu}c^\la_{\mu\nu}
HS_\mu({\bf z};\y)HS_\nu({\bf z}^{-1};{\bf y^{-1}})\Bigl{)}{\rm
ch}V_l^\la.
\end{equation}
On the other hand \thmref{duality} implies that
\begin{equation}\label{character4}
{\rm ch}\F=\sum_{\la}{\rm
ch}L(\hat{gl}_{\infty|\infty},\La(\la)){\rm ch}V_l^\la.
\end{equation}
Using \eqnref{character3} together with \eqnref{character4} we can
prove the following character formula.

\begin{thm}\label{characterthm} We have
\begin{equation*}
{\rm
ch}L(\hat{gl}_{\infty|\infty},\La(\la))=\sum_{\mu,\nu}c^\la_{\mu\nu}
HS_\mu({\bf z};\y)HS_\nu({\bf z}^{-1};{\bf y^{-1}}),
\end{equation*}
where the summation is over all partitions $\mu$ and $\nu$ of
length not exceeding $l$ and $c^\la_{\mu\nu}$ are determined by
the tensor product decomposition
$V^{\mu}_l\otimes(V^{\nu}_l)^*=\sum_{\la}c_{\mu\nu}^\la V^\la_l$.
\end{thm}

\begin{proof}
The statement of the theorem would follow from the linear
independence of the Schur functions in the ring of symmetric
functions, if the summation in \eqnref{character3} and
\eqnref{character4} were over partitons $\la$ of length not
exceeding $l$. However, here we need to deal with summation over
generalized partitions $\la$.  This will follow from the following
lemma.
\end{proof}

\begin{lem}\label{independence}
Let $q$ be an indeterminate and suppose that $\sum_\la\phi(q){\rm
ch}V^\la_l=0$, where $\phi_\la(q)$ are power series in $q$ and
$\la$ above is summed over all generalized partitions of length
not exceeding $l$.  Then $\phi_\la(q)=0$, for all $\la$.
\end{lem}
\begin{proof}
We continue to use the notation from above. We let
$p_{l}=(x_1x_2\cdots x_l)^{\frac{1}{l}}$ and
$p_1=(x_1x_2^{-1})^{\frac{1}{l}},p_2=(x_2x_3^{-1})^{\frac{1}{l}},\cdots,p_{l-1}=
(x_{l-1}x_{l}^{-1})^{\frac{1}{l}}$. Then ${\rm ch}V^\la_l$ may be
written as
\begin{equation}\label{aux11}
{\rm
ch}V^\la_l=p_l^{\sum_{i=1}^l\la_i}\rho_{\la}(p_1,\cdots,p_{l-1}),
\end{equation}
where $\rho_\la(p_1,\cdots,p_{l-1})$ are Laurent polynomials in
$p_1,\cdots,p_{l-1}$. The Laurent polynomial
$\rho_\la(p_1,\cdots,p_{l-1})$ here is of course just the
corresponding character of the irreducible $sl_n$-module.

We need to show that if we have $\sum_{\la}\phi_\la(q){\rm ch
}V^{\la}_l=0$, then $\phi_\la(q)=0$ for all $\la$, where $q$ is
some indeterminate.  Using \eqnref{aux11} by considering just the
coefficient of $p_l^m$, $m\in\Z$, we deduce that
\begin{equation}\label{aux12}
\sum_{\la,\sum\la_i=m}\phi_\la(q)\rho_\la(p_1,\cdots,p_{l-1})=0.
\end{equation}
Now it is clear that $\sum_{i=1}^l\la_i$ and
$\rho_\la(p_1,\cdots,p_{l-1})$ uniquely determines $\la$.  Hence
if we sum over generalized partitions $\la$ with
$\sum_{i=1}^l\la_i$ fixed, $\rho_\la$ is summed over inequivalent
irreducible finite-dimensional $sl_n$-characters.  Thus by the
Weyl character formula we may write
$$\rho_\la(p_1,\cdots,p_{l-1})=\frac{\sum_{w\in
W}\epsilon_we^{w(\la'+\rho)}}{\prod_{\alpha\in\Delta_+}(e^{\alpha/2}-e^{-\alpha/2})},$$
where $\la'$ is the corresponding $sl_l$-highest weight of $\la$,
$W$ is the Weyl group, $\epsilon_w$ is the sign of $w\in W$,
$\Delta_+$ is a set of positive roots and
$\rho=\hf\sum_{\alpha\in\Delta_+}\alpha$.

Multiplying \eqnref{aux12} by
$\prod_{\alpha\in\Delta_+}(e^{\alpha/2}-e^{-\alpha/2})$ and using
the Weyl character formula we get
\begin{equation}\label{aux13}
\sum_{\la'}\phi_\la(q)\sum_{w\in W}e^{w(\la'+\rho)}=0.
\end{equation}
As $\la'+\rho$ is a regular dominant weight, the coefficient of
$e^{\la'+\rho}$ in \eqnref{aux13} above is $\phi_\la(q)$.  Thus
$\phi_\la(q)=0$.
\end{proof}

We conclude this section by applying the character formula to the
case of $l=1$, that is when the central charge is $1$. In this
case $\mu$ and $\nu$ are integers. Furthermore,
$\La(\la)=\la\omega_{\hf}+\La_0$, for $\la\ge 0$, and
$\La(\la)=-\omega_0+(\la+1)\omega_{-\hf}+\La_0$, for $\la<0$.
Since obviously $V^\mu_1\otimes (V^{\nu}_1)^*=V^{\mu-\nu}_1$, we
see that
\begin{align}\label{characterl=1}
&{\rm
ch}L(\hat{gl}_{\infty|\infty},\la\omega_{\hf}+\La_0)=\sum_{\mu-\nu=\la}
HS_\mu({\bf z};\y)HS_\nu({\bf z}^{-1};{\bf y^{-1}}),\quad\la\ge
0,\\
&{\rm
ch}L(\hat{gl}_{\infty|\infty},-\omega_0+(\la+1)\omega_{-\hf}+\La_0)=
\sum_{\mu-\nu=\la} HS_\mu({\bf z};\y)HS_\nu({\bf z}^{-1};{\bf
y^{-1}}),\quad\la< 0,\nonumber
\end{align}
where the summation is over all partitions of non-negative
integers $\mu$ and $\nu$.

\section{A $q$-Character formula for $\hgltwo$-modules}

We can obtain a $q$-character formula for the irreducible highest
weight module $L(\hgltwo,\Lambda(\la))$ from
\thmref{characterthm}. The resulting character formula will be
involve a sum of products of hook Schur functions. However we can
use \propref{combid} and \lemref{independence} to obtain a simpler
formula that will only involve a sum of hook Schur functions. This
we will discuss below.

Since $\hgltwo$ has a principal $\half\Z$-gradation, by declaring
the highest weight vectors to be of degree $0$, its irreducible
highest weight modules are naturally $\half\Z$-graded. Given a
quasi-finite highest weight module $V=\oplus_{s\in\half\Z_+}V_s$,
we can define the {\em $q$-character} of $V$ to be
\begin{equation*}
{\rm ch}_{q}V=\sum_{s\in\half\Z_+}{\rm dim}V_sq^s.
\end{equation*}
It is known \cite{CW3} that $L(\hgltwo,\La(\la))$ is quasi-finite,
so that we may compute its $q$-character.

We first introduce the Virasoro field:
\begin{align}\label{Vir}
L(z)=\sum_{n\in\Z}L_nz^{-n-2}=&\half\sum_{i=1}^n
(:\partial\psi^{+,i}(z)\psi^{-,i}(z):-
:\psi^{+,i}(z)\partial\psi^{-,i}(z):)\\
+&\half\sum_{i=1}^n(:\gamma^{+,i}(z)\partial\gamma^{-,i}(z):-
:\partial\gamma^{+,i}(z)\gamma^{-,i}(z):).\nonumber
\end{align}

Now set $\tilde{d}=-L_0-\half\alpha_0$, where $\alpha_0=
\sum_{i=1}^l(\sum_{n\in\Z}:\psi^{+,i}_{-n}\psi^{-,i}_{n}:)$. We
have
\begin{equation*}
[\tilde{d},\psi^{\pm,i}_n]=n\psi^{\pm,i}_n,\quad
[\tilde{d},\gamma^{\pm,i}_r]=r\gamma^{\pm,i}_r.
\end{equation*}
It is clear that $\tilde{d}$ commutes with $\sum_{i=1}^l\C E_{ii}$
so that we may decompose the Fock space $\F$ into its
$\C\tilde{d}+\sum_{i=1}^l\C E_{ii}$-weight spaces. Let ${\rm
ch}_q\F$ denote the resulting character. Letting
$x_i=e^{\epsilon_i}$ be as before and $e^{-\delta'}=q$ with
$\delta'$ defined by $\delta'(E_{ii})=0$ and
$\delta'(\tilde{d})=1$ we have
\begin{equation*}
{\rm ch}_{q}\F=\prod_{i=1}^l\Bigl{(}\prod_{n\in\N,r\in\half+\Z_+}
\frac{(1+q^nx_i)(1+q^nx_i^{-1})}{(1-q^rx_i)(1-q^rx_i^{-1})}(1+x_i^{-1})\Bigl{)}.
\end{equation*}

By \propref{combid} we have
\begin{equation}\label{aux17}
\prod_{i=1}^l\prod_{n\in\N,r\in\half+\Z_+}
\frac{(1+q^nx_i)(1+q^nx_i^{-1})}{(1-q^rx_i)(1-q^rx_i^{-1})}=
\sum_{\mu}HS_\mu({\bf q}^r;{\bf q}^n)s_\mu(\x,\x^{-1}),
\end{equation}
where ${\bf q}^r=\{q^{\half},q^{\frac{3}{2}},\cdots\}$, ${\bf
q}^n= \{q,q^2,\cdots\}$.  Here the summation of $\mu$ is over all
partitions of length not exceeding $2l$, and
$s_\mu(\x,\x^{-1})=s_{\mu}(x_1,\cdots,x_l,x_1^{-1},\cdots,x_l^{-1})$.
The expression $s_\mu(\x,\x^{-1})$ has a simple interpretation.
Consider the embedding of the Lie algebra $gl_l$ into $gl_{2l}$
given by
\begin{equation*}
A\in gl_l\rightarrow
\begin{pmatrix}
A&0\\
0&-A^{t}
\end{pmatrix}\in gl_{2l}.
\end{equation*}
The $gl_{2l}$-module $V^\mu_{2l}$ has $gl_{2l}$-character
$s_\mu(x_1,x_2,\cdots,x_{2l})$.  We may restrict $V^\mu_{2l}$ via
the embedding above to a $gl_l$-module.  Then its $gl_l$-character
is given by $s_\mu(\x,\x^{-1})$.

Of course the product $\prod_{i=1}^l (1+x_i^{-1})$ is nothing but
the character of $\Lambda^{\bullet}(\C^{l*})$, the exterior
algebra of the $gl_l$-module contragredient to the standard
module. We may decompose
$V^{\mu}_{2l}\otimes\Lambda^{\bullet}(C^{l*})$ as a $gl_l$-module
and obtain
\begin{equation*}
V^{\mu}_{2l}\otimes\Lambda^{\bullet}(C^{l*})\cong\sum_{\mu}c_\la^\mu
V^\la_l,
\end{equation*}
where $\la$ are now generalized partitions of lengths not
exceeding $l$ and $c^\mu_\la$ are the multiplicities of this
tensor product decomposition. Therefore we have the corresponding
character identity
\begin{equation}\label{aux18}
s_\mu(\x,\x^{-1})\prod_{i=1}^l(1+x_i^{-1})=\sum_{\la}c_{\la}^\mu
{\rm ch}V^\la_l.
\end{equation}
Using \eqnref{aux17} and \eqnref{aux18} we obtain
\begin{equation}\label{aux19}
{\rm ch}_q\F=\sum_{\la}\Bigl{(}\sum_\mu c^\mu_\la HS_\mu({\bf
q}^r;{\bf q}^n)\Bigl{)}{\rm ch}V^\la_l,
\end{equation}
where $\mu$ here is summed over all generalized partitions of
length not exceeding $l$.

On the other hand by \thmref{duality}, using the explicit formula
of the joint highest weight vectors, we see that
\begin{equation}\label{aux20}
{\rm ch}_q\F=\sum_{\la}q^{h(\la)}{\rm ch}_q
L(\hgltwo,\La(\la)){\rm ch}V^\la_l,
\end{equation}
where $h(\la)=\sum_{r\in\hf\Z}r\La(\la)_r$.  Combining
\eqnref{aux19} with \eqnref{aux20}, using \lemref{independence},
we have thus proved the following.
\begin{thm}
\begin{equation*}
{\rm ch}_q L(\hgltwo,\La(\la))=q^{-h(\la)}\sum_{\mu}c^\la_\mu
HS_\la({\bf q}^r,{\bf q}^n),
\end{equation*}
where the sum is over all partitions of length not exceeding $2l$
and the coefficient $c^\mu_\la$ are determined by \eqnref{aux18}.
\end{thm}

We consider the simplest case when $l=1$ and let $x=x_1$. In this
case it is clear that $s_\mu(x,x^{-1})$ is just the character of
the corresponding $sl_2$-module, and hence for a partitions
$\mu=(\mu_1,\mu_2)$ with $\mu_1\ge\mu_2\ge 0$ we have
\begin{equation*}
s_\mu(x,x^{-1})=x^{\mu_1-\mu_2}+x^{\mu_1-\mu_2-2}+\cdots+x^{\mu_2-\mu_1}.
\end{equation*}
Hence
\begin{equation*}
s_\mu(x,x^{-1})(1+x^{-1})=x^{\mu_1-\mu_2}+x^{\mu_1-\mu_2-1}+\cdots+x^{\mu_2-\mu_1-1}.
\end{equation*}
Therefore we see that for $\la\in\Z$
\begin{align*}
&{\rm
ch}_qL(\hgltwo,\la\omega_{\hf}+\La_0)=q^{-\frac{\la}{2}}\sum_{\mu_1-\mu_2\ge\la}HS_{\mu}({\bf
q}^r;{\bf q}^n),\quad\la\ge 0,\\
&{\rm ch}_qL(\hgltwo,-\omega_0+(\la+1)\omega_{-\hf}+\La_0)=
q^{\frac{\la+1}{2}}\sum_{\mu_2-\mu_1-1\le\la}HS_{\mu}({\bf
q}^r;{\bf q}^n),\quad\la< 0.
\end{align*}
Note that the coefficient of $q^s$, $s\in\half\Z_+$ in
$HS_{\mu}({\bf q}^r;{\bf q}^n)$ can be computed as follows.
Arrange ${\bf q}^r=\{q^{\half}<q^{\frac{3}{2}}<\cdots\}$ and ${\bf
q}^n=\{q<q^2<\cdots\}$ in increasing order. Let $T$ be an
$(\infty|\infty)$-semi-standard tableau of shape $\mu$. Let
$q^{m(T)}$ denote the product of all entries in $T$ so that
$m(T)\in\half\Z_+$. Then for a fixed $s\in\half\Z_+$ the
coefficient of $q^s$ in $HS_\mu({\bf q}^r;{\bf q}^n)$ is the
number of $(\infty|\infty)$-semi-standard tableaux of shape $\mu$
with $m(T)=s$.  Hence we have the formula
\begin{equation*}
HS_{\mu}({\bf q}^r;{\bf
q}^n)=\sum_{s\in\half\Z_+}(\sum_{m(T)=s}1)q^s.
\end{equation*}
For example if $\mu$ is the partition $(2,0)$, then
$HS_{(2,0)}({\bf q}^r,{\bf q}^n)=q+q^{\frac{3}{2}}+2q^2+2
q^{\frac{5}{2}}+\cdots$.

Now by \cite{KL}
\begin{align*}
&{\rm ch}_qL(\hgltwo,\la\omega_{\hf}+\La_0)=
(1+q^{\la+\half})^{-1}\prod_{n\in\N,r\in\half+\Z_+}
\frac{(1+q^r)^2}{(1-q^n)^2},\quad
\la\ge 0,\\
&{\rm ch}_qL(\hgltwo,-\omega_0+(\la+1)\omega_{-\hf}+\La_0)=
(1+q^{-\la-\half})^{-1}\prod_{n\in\N,r\in\half+\Z_+}
\frac{(1+q^r)^2}{(1-q^n)^2},\quad\la<0.
\end{align*}

Therefore we obtain the following interesting combinatorial
identities.

\begin{thm} For $\la\in\Z$ and $\mu=(\mu_1,\mu_2)$ a partition of length at most
$2$ we have
\begin{align*}
&\sum_{\mu_1-\mu_2\ge\la}HS_{\mu}({\bf q}^r;{\bf q}^n) =
q^{\frac{\la}{2}}(1+q^{\la+\half})^{-1}\prod_{n\in\N,r\in\half+\Z_+}
\frac{(1+q^r)^2}{(1-q^n)^2},\quad\la\ge 0,\\
&\sum_{\mu_2-\mu_1-1\le\la}HS_{\mu}({\bf q}^r;{\bf q}^n) =
q^{-\frac{\la+1}{2}}(1+q^{-\la-\half})^{-1}\prod_{n\in\N,r\in\half+\Z_+}
\frac{(1+q^r)^2}{(1-q^n)^2},\quad\la<0.
\end{align*}
\end{thm}

First we note that ${\rm
ch}_qL(\hgltwo,\La_0)=\sum_{\mu}HS_{\mu}({\bf q}^r;{\bf q}^n)$,
where the summation is over all $\mu$.  Also we have
\begin{equation}\label{qlevel1}
{\rm ch}_q L(\hgltwo,\la\omega_{\hf}+\La_0)={\rm ch}_q L
(\hgltwo,-\omega_0-\la\omega_{-\hf}+\La_0),
\end{equation}
for all $\la\in\Z_+$.  We have a similar $q$-character identity
for higher level representations as well. In fact one has the
following proposition.
\begin{prop}  Let $l\in\N$ and $\la$ be a generalized partition of length not
exceeding $l$.  We have
\begin{align*}
{\rm ch}_q L(\hgltwo,\La(\la_1,\la_2,\cdots,\la_l)&)=\\
&{\rm ch}_q L (\hgltwo,\La(-\la_l-1,\cdots,-\la_2-1,-\la_1-1)).
\end{align*}
In particular when $l=1$ the above identity reduces to
\eqnref{qlevel1}.
\end{prop}
\begin{proof}
We will use the Howe duality involving a $B$-type subalgebra of
$\hgltwo$ and the double covering Lie group ${\rm Pin}_{2l}$ of
the spin group ${\rm Spin}_{2l}$ on $\F$ in \cite{CW3} to prove
this identity.

We recall that the subalgebra $\ospd$ is the subalgebra of
$\gltwo$ preserving the even non-degenerate bilinear form
$(\cdot|\cdot)$ of $\gltwo$ defined by
\begin{align*}
&(e_{i}|e_j)=(-1)^{i}\delta_{i,-j}, \quad i,j\in\Z,\\
&(e_{r}|e_s)=(-1)^{r+\hf}\delta_{r,-s},\quad r,s\in\hf+\Z.
\end{align*}
Restricting the $2$-cocycle of $\gltwo$ to $\ospd$ we obtain a
central extension $\hospd$ of $\ospd$. Obviously $\hospd$ acts on
$\F$ and in fact there exists an action of ${\rm Pin}_{2l}$ on
$\F$ such that $(\hospd,{\rm Pin}_{2l})$ forms a dual pair
\cite{CW3}. Hence we have a multiplicity-free decomposition
\begin{equation*}
\F\cong \sum_{\mu}L(\hospd,\mu')\otimes W_{2l}^\mu,
\end{equation*}
where $W_{2l}^\mu$ ( respectively $L(\hospd,\mu')$) stands for an
irreducible ${\rm Pin}_{2l}$-module (respectively $\hospd$-module)
of highest weight $\mu$ (respectively $\mu'$) and the map
$\mu\rightarrow\mu'$ is a bijection. On the other hand by
\thmref{duality} we have another dual pair $(\hgltwo\times GL_l)$
on $\F$.  Now if we twist the action of $GL_l$ on $\F$ by $({\rm
det})^{\hf}$, then these two dual pairs form a seesaw pair in the
sense of Kudla \cite{Ku}.  This implies that we have the following
decompositions of respectively $\hospd$- and $GL_l$-modules:
\begin{align*}
L(\hgltwo,\La(\la'))&\cong\bigoplus_{\mu}b_\la^\mu L(\hospd,\mu'),\\
W_{2l}^\mu&\cong\bigoplus_{\la}b_\la^\mu V_l^\la,
\end{align*}
where here $\la'=\la-(\hf,\cdots,\hf)$ due to the twist by ${\rm
det}^\hf$ and $b_\la^\mu\in\Z_+$.

Now the ${\rm Pin}_{2l}$-modules that appear in the decomposition
of $\F$ with respect to the dual pair $(\hospd,{\rm Pin}_{2l})$,
when regarded as a module over the Lie algebra $so_{2l}$,
decomposes into two irreducible modules contragredient to each
other.  Hence $b_{\la}^{\mu}=b_{\la^*}^{\mu}$.  But then as
$\hospd$-modules we have
\begin{align*}
L(\hgltwo,\La(\la'))&\cong\bigoplus_{\mu}b_\la^\mu L(\hospd,\mu')\\
&\cong\bigoplus_{\mu}b_{\la^*}^\mu L(\hospd,\mu')\\&\cong
L(\hgltwo,\La({\la^*}')).
\end{align*}
Hence as modules over $\hospd$ we have an isomorphism
\begin{equation*}
L(\hgltwo,\La(\la_1,\la_2,\cdots,\la_l))\cong L
(\hgltwo,\La(-\la_l-1,\cdots,-\la_2-1,-\la_1-1)).
\end{equation*}
Since $\hospd$ is a subalgebra of $\hgltwo$ preserving the
principal $\Z$-gradation, the proposition follows.
\end{proof}

\section{A Character formula for $\hat{gl}_{m|n}$-modules at level $1$ }

We first recall a method of constructing representations of affine
superalgebras $\hat{gl}_{n|n}$ from representations of
$\hat{gl}_{\infty|\infty}$ \cite{KL}. It is a generalization of
the classical reduction from $\hat{gl}_{\infty}$ to $\hat{gl}_n$.
Our presentation below is somewhat different from \cite{KL} in
flavor.

Let $\psi^{+}(z)=\sum_{n\in\Z}\psi^+_nz^{-n-1}$ and
$\psi^-(z)=\sum_{n\in\Z}\psi^-_nz^{-n}$ be a pair of free fermions
and let
$\gamma^{\pm}(z)=\sum_{r\in\frac{1}{2}+\Z}\gamma^{\pm}_rz^{-r-\frac{1}{2}}$
be a pair of free bosons. Let $\F$ denote the corresponding Fock
space generated by the vacuum vector $|0>$.  We have thus an
action of the dual pair $(\hat{gl}_{\infty|\infty},gl_1)$ on $\F$,
where the central charge of $\hat{gl}_{\infty|\infty}$ is $1$.
From the pair of free fermions we may construct $n$ pairs of free
fermions $\psi^{\pm,i}(z)$, for $i=1,\cdots,n$, as follows:
\begin{align}
\psi^{+,i}(z)&=\sum_{k\in\Z}\psi^{+,i}_kz^{-k-1}=\sum_{k\in\Z}\psi^+_{-i+n(k+1)}z^{-k-1},\\
\psi^{-,i}(z)&=\sum_{k\in\Z}\psi^{-,i}_kz^{-k}=\sum_{k\in\Z}\psi^-_{i+n(k-1)}z^{-k}.
\end{align}
It is easy to check that the only non-zero commutation relations
are
\begin{equation*}
[\psi^{+,i}_m,\psi^{-,j}_n]=\delta_{ij}\delta_{m+n,0},\quad
m,n\in\Z.
\end{equation*}
Similarly we construct from our pair of free bosons $n$ pairs of
free bosons $\gamma^{\pm,i}(z)$ for $i=1,\cdots,n$, via
\begin{align}
\gamma^{+,i}(z)&=\sum_{r\in\frac{1}{2}+\Z}\gamma^{+,i}_rz^{-r-\frac{1}{2}}=
\sum_{r\in\frac{1}{2}+\Z}\gamma^+_{-i+\frac{1}{2}+n(r+\frac{1}{2})}z^{-r-\half},\\
\gamma^{-,i}(z)&=\sum_{r\in\frac{1}{2}+\Z}\gamma^{-,i}_rz^{-r-\frac{1}{2}}=
\sum_{r\in\frac{1}{2}+\Z}\gamma^-_{i-\frac{1}{2}+n(r-\frac{1}{2})}z^{-r-\half}.
\end{align}
Again it is easily checked that the only non-zero commutation
relations are
\begin{equation*}
[\gamma^{+,i}_r,\gamma^{-,j}_s]=\delta_{ij}\delta_{r+s,0},\quad
r,s\in\frac{1}{2}+\Z.
\end{equation*}

We may now use these $n$ pairs of fermions and bosons to construct
a copy of the affine $gl_{n|n}$ of central charge $1$ in the
standard way:

\begin{align*}
E_{ij}(z)&=\sum_{m\in\Z}E_{ij}(m)z^{-m-1}=:\psi^{+,i}(z)\psi^{-,j}(z):,\\
e_{ij}(z)&=\sum_{m\in\Z}e_{ij}(m)z^{-m-1}=-:\gamma^{+,i}(z)\gamma^{-,j}(z):,\\
a_{ij}(z)&=\sum_{m\in\Z}a_{ij}(m)z^{-m-1}=:\psi^{+,i}(z)\gamma^{-,j}(z):,\\
b_{ij}(z)&=\sum_{m\in\Z}b_{ij}(m)z^{-m-1}=-:\gamma^{+,i}(z)\psi^{-,j}(z):.
\end{align*}
Explicitly we have the following formulas.
\begin{align*}
E_{ij}(m)&=\sum_{k+l=m}:\psi^+_{-i+n(k+1)}\psi^-_{j+n(l-1)}:,\\
e_{ij}(m)&=-\sum_{r+s=m}:\gamma^+_{-i+\frac{1}{2}+n(r+\half)}\gamma^-_{j-\half+n(s-\half)}:,\\
a_{ij}(m)&=\sum_{k+r+\half=m}:\psi^+_{-i+n(k+1)}\gamma^-_{j-\half+n(r-\half)}:,\\
b_{ij}(m)&=-\sum_{k+r-\half=m}:\gamma^+_{\half-i+n(r+\half)}\psi^-_{j+n(k-1)}:.
\end{align*}

The following lemma is straightforward.

\begin{lem}
We have
\begin{align*}
&E_{ij}(m)\in (\hat{gl}_{\infty|\infty})_+,\ {\rm for}\ i<j {\rm\
with\ } m=0 {\rm\ and\ for\ }
m\ge 1 ,\\
&e_{ij}(m)\in (\hat{gl}_{\infty|\infty})_+,\ {\rm for}\ i<j {\rm\
with\ } m=0 {\rm\ and\ for\ }
m\ge 1,\\
&a_{ij}(m)\in (\hat{gl}_{\infty|\infty})_+,\ {\rm for}\ i<j {\rm\
with\ } m=0 {\rm\ and\ for\ }
m\ge 1,\\
&b_{ij}(m)\in (\hat{gl}_{\infty|\infty})_+,\ {\rm for}\ i\le j
{\rm\ with\ } m=0 {\rm\ and\ for\ }m\ge 1.
\end{align*}
\end{lem}

Now every representation of $\hat{gl}_{\infty|\infty}$ that
appears in the decomposition of $\F$ is obviously invariant under
the action of $\hat{gl}_{n|n}$ constructed via reduction modulo
$n$. By \cite{KL} every irreducible
$\hat{gl}_{\infty|\infty}$-module that appears in $\F$ in fact
remains irreducible when restricted to $\hat{gl}_{n|n}$. Hence it
follows from the previous lemma that every irreducible
representation of $\hat{gl}_{\infty|\infty}$ that appears in the
decomposition of $\F$ is a highest weight irreducible
$\hat{gl}_{n|n}$-module with respect to the Borel subalgebra
induced by the non-standard Borel subalgebra of $gl_{n|n}$.

\begin{rem}
In general one can construct $m$ pairs of free fermions and $n$
pairs of free bosons using the method just described. One then
constructs a copy of the affine $gl_{m|n}$ of central charge $1$
in the usual way. However, the resulting representations of this
affine $gl_{m|n}$ on the highest weight irreducible
representations of $\hat{gl}_{\infty|\infty}$ of are not highest
weight representations with respect to a Borel subalgebra induced
from a Borel subalgebra of $gl_{m|n}$.
\end{rem}

We want to deduce a character formula for the
$\hat{gl}_{n|n}$-module $L(\hat{gl}_{\infty|\infty},\La(\la))$
from \eqnref{characterl=1}. For this we will need to find a
slightly more general formula for ${\rm
ch}L(\hat{gl}_{\infty|\infty},\La(\la))$.

Set $\bar{d}=-L_0$ (see \eqnref{Vir}) so that we have
\begin{equation*}
[\bar{d},\psi^{\pm,i}_k]=(k\pm\half)\psi^{\pm,i}_k,\quad
[\bar{d},\gamma^{\pm,i}_r]=r\gamma^{\pm,i}_r.
\end{equation*}
From this the following lemma is a straightforward computation.
\begin{lem}\label{aux14}
For all $m\in\Z$ we have $[\bar{d},X(m)]=mX(m)$, for $X=E,e,a,b$.
\end{lem}

By construction we see that $\bar{d}$ acts on each
$L(\hat{gl}_{\infty|\infty},\La(\la))$ as a semisimple linear
operator and also $[\bar{d},(\hat{gl}_{\infty|\infty})_0]=0$.
Hence we may compute the character of
$L(\hat{gl}_{\infty|\infty},\La(\la))$ with respect to
$(\hat{gl}_{\infty|\infty})_0\oplus\C\bar{d}$.  Letting $\delta\in
((\hat{gl}_{\infty|\infty})_0\oplus\C\bar{d})^*$ defined by
$\delta(\bar{d})=1$ and $\delta((\hat{gl}_{\infty|\infty})_0)=0$
we have ($x=x_1$, $q=e^{-\delta}$):
\begin{equation*}
{\rm ch}\F=\prod_{k\in\N,s\in\half+\Z_+}\prod_{j=1}^n\frac{(1+x
y_{j+n(k-1)}q^{k-\half}) (1+x^{-1}y_{j-nk}^{-1}q^{k-\half})}{(1-x
z_{j-\half+(s-\half)n}q^s)(1-x^{-1}z_{-j+\half-(s-\half)n}^{-1}q^s)}.
\end{equation*}
It follows again from \propref{combid} that with respect to
$(\hgltwo)_0\oplus\C d$ the character of
$L(\hat{gl}_{\infty|\infty},\La(\la))$ equals to
\begin{equation}\label{aux15}
\sum_{\mu-\nu=\la} HS_\mu({\bf zq};\y{\bf q})HS_\nu({\bf
z}^{-1}{\bf q};{\bf y^{-1}q}),
\end{equation}
where for $k\in\N$, $s\in\half+\Z_+$ and $j=1,\cdots,n$. Here
$\y{\bf q}=\{y_{j+n(k-1)}q^{k-\half}\}$, $\y^{-1}{\bf
q}=\{y^{-1}_{j-nk}q^{k-\half}\}$, ${\bf
zq}=\{z_{j-\half+(s-\half)n}q^s\}$ and ${\bf
z^{-1}q}=\{z^{-1}_{{-j+\half}-(s-\half)n}q^s\}$.

The following lemma follows easily from our construction.

\begin{lem}\label{aux16} For $i,j=1,\cdots,n$; $k\in\Z$ and $r\in\half+\Z$ we have
\begin{align*}
&[E_{ii}(0),\psi^{\pm,j}_k]=\pm\delta_{ij}\psi^{\pm,j}_k,\quad
[E_{ii}(0),\gamma^{\pm,j}_r]=0,\\
&[e_{ii}(0),\psi^{\pm,j}_k]=0,\quad
[e_{ii}(0),\gamma^{\pm,j}_r]=\pm\delta_{ij}\gamma^{\pm,j}_{r}.
\end{align*}
\end{lem}

Let us denote the scaling operator of $\hat{gl}_{n|n}$ by $d$.  It
is clear from \lemref{aux14} that the linear operator $d-\bar{d}$
acts as a scalar on each irreducible representation of
$\hat{gl}_{n|n}$ that appears in the decomposition of $\F$.  The
scalar can be computed from the explicit formulas of the
$\hat{gl}_{\infty|\infty}$ highest weight vectors given by
\thmref{duality}.  In the case of $l=1$ the highest weight vectors
are given by
\begin{align}
&(\gamma^+_{-\half})^{\la}|0>,\quad\la\ge 0,\label{positivela}\\
&(\gamma^-_{-\half})^{-\la-1}\psi^-_0|0>,\quad\la<0\label{negativela}.
\end{align}
Now $\bar{d}$ acts on the former with eigenvalue $-\frac{\la}{2}$
and on the latter with eigenvalue $\frac{\la}{2}$. Thus using
\lemref{aux16} we can rewrite \eqnref{aux15} in the following
form.

\begin{thm}\label{main1} For $\y{\bf q}=\{y_j q^s\}$, $\y^{-1}{\bf q}=\{y^{-1}_j
q^s\}$, ${\bf zq}=\{z_jq^s\}$ and ${\bf z^{-1}q}=\{z^{-1}_jq^s\}$
with $j=1,\cdots,n$ and $s\in\half+\Z_+$ we have ($\mu,\nu\in\Z_+$
and $\la\in\Z$)
\begin{align*}
&{\rm ch}L(\hat{gl}_{n|n},\la\tilde{\delta}_1+\tilde{\Lambda}_0)=
q^{-\frac{\la}{2}}\sum_{\mu-\nu=\la} HS_\mu({\bf zq};\y{\bf
q})HS_\nu({\bf z}^{-1}{\bf q};{\bf
y^{-1}q}),\quad\la\ge 0,\\
&{\rm
ch}L(\hat{gl}_{n|n},(\la+1)\tilde{\delta}_{n}-\tilde{\epsilon}_n+\tilde{\Lambda}_0)=
q^{\frac{\la}{2}}\sum_{\mu-\nu=\la}HS_\mu({\bf zq};\y{\bf
q})HS_\nu({\bf z}^{-1}{\bf q};{\bf y^{-1}q}),\quad\la< 0.
\end{align*}
where $y_j=e^{\tilde{\epsilon}_j}$, $z_j=e^{\tilde{\delta}_j}$ and
$q=e^{-\tilde{\delta}}$.
\end{thm}

\begin{prop}
The $\hat{gl}_{n|n}$-modules in \thmref{main1} are integrable.
Furthermore  they form a complete list of integrable
$\hat{gl}_{n|n}$-modules of level $1$.
\end{prop}

\begin{proof} We will employ the method of simple odd
reflections \cite{PS} following \cite{KW2}.

Recall that the set of simple roots of $\hat{gl}_{n|n}$ with
respect to the standard Borel subalgebra is given by
$\{\alpha_0=\tilde{\delta}_n-\tilde{\epsilon}_1+\tilde{\delta},
\alpha_1=\tilde{\epsilon}_1-\tilde{\epsilon}_2,
\alpha_2=\tilde{\epsilon}_2-\tilde{\epsilon}_3,\cdots,
\alpha_n=\tilde{\epsilon}_n-\tilde{\delta}_1,
\alpha_{n+1}=\tilde{\delta}_1-\tilde{\delta}_2, \cdots,
\alpha_{2n-1}=\tilde{\delta}_{n-1}-\tilde{\delta}_n\}$.  The
corresponding Dynkin diagram is as follows (as usual $\bigotimes$
denotes an isotropic odd root):
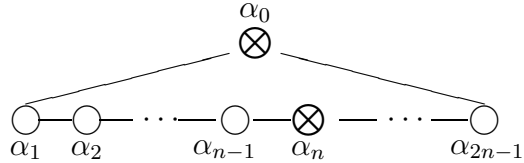
\begin{table}[hb]
\vspace*{-8ex}$
\begin{array}{c c}
\setlength{\unitlength}{0.16in}
\begin{picture}(20,8)
\put(0,1.5){\line(4,1){6.65}} \put(15,1.5){\line(-4,1){6.65}}
\put(7.5,4.5){\makebox(0,0)[c]{$\alpha_0$}}
\put(7.5,3.5){\makebox(0,0)[c]{$\bigotimes$}}
\put(0,1){\makebox(0,0)[c]{$\bigcirc$}}
\put(2,1){\makebox(0,0)[c]{$\bigcirc$}}
\put(6.85,1){\makebox(0,0)[c]{$\bigcirc$}}
\put(9.25,1){\makebox(0,0)[c]{$\bigotimes$}}
\put(15,1){\makebox(0,0)[c]{$\bigcirc$}}
\put(0.4,1){\line(1,0){1.1}} \put(2.4,1){\line(1,0){1.1}}
\put(5,1){\line(1,0){1.2}} \put(7.45,1){\line(1,0){1.2}}
\put(10.25,1){\line(1,0){1.2}} \put(13.25,1){\line(1,0){1.2}}
\put(4.5,1){\makebox(0,0)[c]{$\cdots$}}
\put(12.5,1){\makebox(0,0)[c]{$\cdots$}}
\put(0,0){\makebox(0,0)[c]{$\alpha_1$}}
\put(2,0){\makebox(0,0)[c]{$\alpha_2$}}
\put(6.5,0){\makebox(0,0)[c]{$\alpha_{n-1}$}}
\put(9.25,0){\makebox(0,0)[c]{$\alpha_n$}}
\put(15.1,0){\makebox(0,0)[c]{$\alpha_{2n-1}$}}
\end{picture}
\end{array}$
\end{table}

Now the set of simple roots with respect to the non-standard Borel
subalgebra are given by
$\{\beta_0=\tilde{\epsilon}_n-\tilde{\delta}_1+\tilde{\delta},
\beta_1=\tilde{\delta}_1-\tilde{\epsilon}_1,
\beta_2=\tilde{\epsilon}_1-\tilde{\delta}_2,
\beta_3=\tilde{\delta}_2-\tilde{\epsilon}_2,\cdots,
\beta_{2n-2}=\tilde{\epsilon}_{n-1}-\tilde{\delta}_n,
\beta_{2n-1}=\tilde{\delta}_{n}-\tilde{\epsilon}_n\}$, with the
corresponding Dynkin diagram
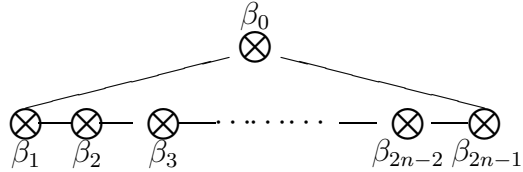
\begin{table}[hb]
\vspace*{-8ex}$
\begin{array}{c c}
\setlength{\unitlength}{0.16in}
\begin{picture}(20,8)
\put(0,1.5){\line(4,1){6.65}} \put(15,1.5){\line(-4,1){6.65}}
\put(7.5,4.5){\makebox(0,0)[c]{$\beta_0$}}
\put(7.5,3.5){\makebox(0,0)[c]{$\bigotimes$}}
\put(0,1){\makebox(0,0)[c]{$\bigotimes$}}
\put(2,1){\makebox(0,0)[c]{$\bigotimes$}}
\put(6.85,1){\makebox(0,0)[c]{$\cdots$}}
\put(9.25,1){\makebox(0,0)[c]{$\cdots$}}
\put(15,1){\makebox(0,0)[c]{$\bigotimes$}}
\put(0.4,1){\line(1,0){1.1}} \put(2.4,1){\line(1,0){1.1}}
\put(5,1){\line(1,0){1.2}} \put(8,1){\makebox(0,0)[c]{$\cdots$}}
\put(10.25,1){\line(1,0){1.2}} \put(13.25,1){\line(1,0){1.2}}
\put(4.5,1){\makebox(0,0)[c]{$\bigotimes$}}
\put(12.5,1){\makebox(0,0)[c]{$\bigotimes$}}
\put(12.5,0){\makebox(0,0)[c]{$\beta_{2n-2}$}}
\put(4.5,0){\makebox(0,0)[c]{$\beta_3$}}
\put(0,0){\makebox(0,0)[c]{$\beta_1$}}
\put(2,0){\makebox(0,0)[c]{$\beta_2$}}
\put(15.1,0){\makebox(0,0)[c]{$\beta_{2n-1}$}}
\end{picture}
\end{array}$
\end{table}

We can use a chain of odd reflections \cite{PS} to bring the
second diagram to the first as follows.  We reflect first along
the odd simple root $\beta_{2n-1}$. After that we reflect along
the right most odd simple root in the bottom row of the Dynkin
diagram which is of the form
$\tilde{\delta}_i-\tilde{\epsilon}_j$. For example, after
reflecting along $\beta_{2n-1}$ we obtain the diagram
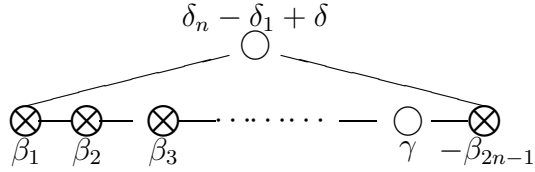
\begin{table}[hb]
\vspace*{-8ex}$
\begin{array}{c c}
\setlength{\unitlength}{0.16in}
\begin{picture}(20,8)
\put(0,1.5){\line(4,1){6.65}} \put(15,1.5){\line(-4,1){6.65}}
\put(7.5,4.5){\makebox(0,0)[c]{$\tilde{\delta}_n-\tilde{\delta}_1+\tilde{\delta}$}}
\put(7.5,3.5){\makebox(0,0)[c]{$\bigcirc$}}
\put(0,1){\makebox(0,0)[c]{$\bigotimes$}}
\put(2,1){\makebox(0,0)[c]{$\bigotimes$}}
\put(6.85,1){\makebox(0,0)[c]{$\cdots$}}
\put(9.25,1){\makebox(0,0)[c]{$\cdots$}}
\put(15,1){\makebox(0,0)[c]{$\bigotimes$}}
\put(0.4,1){\line(1,0){1.1}} \put(2.4,1){\line(1,0){1.1}}
\put(5,1){\line(1,0){1.2}} \put(8,1){\makebox(0,0)[c]{$\cdots$}}
\put(10.25,1){\line(1,0){1.2}} \put(13.25,1){\line(1,0){1.2}}
\put(4.5,1){\makebox(0,0)[c]{$\bigotimes$}}
\put(12.5,1){\makebox(0,0)[c]{$\bigcirc$}}
\put(12.5,0){\makebox(0,0)[c]{$\gamma$}}
\put(4.5,0){\makebox(0,0)[c]{$\beta_3$}}
\put(0,0){\makebox(0,0)[c]{$\beta_1$}}
\put(2,0){\makebox(0,0)[c]{$\beta_2$}}
\put(15.1,0){\makebox(0,0)[c]{$-\beta_{2n-1}$}}
\end{picture}
\end{array}$
\end{table}

\noindent where
$\gamma=\tilde{\epsilon}_{n-1}-\tilde{\epsilon}_n$.  As
$-\beta_{2n-1}$ and $\gamma$ are not of the form
$\tilde{\delta}_i-\tilde{\epsilon}_j$, the next step is to reflect
along the odd simple root $\beta_{2n-3}$.  Continuing this way we
obtain the diagram corresponding to the standard Borel subalgebra.

Now according to Lemma 1.4 of \cite{KW2} a highest weight vector
$v_{\La}$ of highest weight $\La$ with respect to the original
Borel subalgebra remains a highest weight vector of the new Borel
subalgebra if and only if the $\La(\check{\gamma})=0$, where
$\check{\gamma}$ is the simple coroot corresponding to the odd
simple root $\gamma$, along which we have reflected. Furthermore
in this case the new highest weight and the original highest
weight coincide. If however $\La(\check{\gamma})\not=0$, then
$e_{-\gamma}v$ is the highest weight vector with respect to the
new Borel subalgebra, where $e_{-\gamma}$ is the root vector
corresponding to $-\gamma$. Furthermore the new highest weight is
$\La-\gamma$.

Let $0\le \la\le n$ and consider the highest weight
$\la\tilde{\delta}_1+\tilde{\La}_0$.  It follows that when we
change from non-standard to the standard Borel subalgebra it gets
transformed to
$\tilde{\epsilon}_1+\cdots+\tilde{\epsilon}_\la+\tilde{\La}_0$. If
$\la> n$, then the highest weight gets transformed to
$\tilde{\epsilon}_1+\cdots+\tilde{\epsilon}_n+
(\la-n)\tilde{\delta}_1+\tilde{\La}_0$. On the other hand let
$\la<0$ and consider the highest weight of the form
$(\la+1)\tilde{\delta}_n+\tilde{\epsilon}_n+\tilde{\La}_0$.
Changing from the non-standard Borel to the standard Borel via the
sequence of odd reflections described above it follows that the
highest weight gets transformed to
$\la\tilde{\delta}_n+\tilde{\La}_0$. However, the list of highest
weights here for the standard Borel subalgeba coincides with the
list of integrable highest weights in \cite{KW2}.  Thus all our
modules are integrable.
\end{proof}

As the Cartan subalgebras of the non-standard Borel and standard
Borel subalgebras coincide, it follows that our character formula
agree with the character formula of \cite{KW2}.  Comparing both
formulas gives rise to combinatorial identities.

Our method can also be used to obtain a character formula for
integrable level $1$ $\hat{gl}_{m|n}$-modules as follows. Consider
the Fock space $\F$ generated by $m$ pairs of free fermions
$\psi^{\pm,i}(z)$, $i=1,\cdots,m$ and $n$ pairs of free bosons
$\gamma^{\pm,j}(z)$, $j=1,\cdots,n$. Then $\F$ according to
\cite{KW2} is completely reducible as a level $1$
$\hat{gl}_{m|n}$-module. To be more precise there is an action of
$\hat{gl}_{m|n}\times gl_1$ on $\F$ and they form a dual pair. The
action of the $\hat{gl}_{m|n}$ is given in the usual way, while
$gl_1$ is generated by the {\em charge operator}
$I=\sum_{i=1}^m\Bigl{(}\sum_{s\in\hf+\Z_+}:\psi^{+,i}_{-s}\psi^{-,i}_{s}:\Bigl{)}-
\sum_{j=1}^n\Bigl{(}\sum_{r\in\hf+\Z_+}:\gamma^{+,j}_{-r}\gamma^{-,j}_{r}:\Bigl{)}$.
We have the following decomposition with respect to this joint
action \cite{KW2}:
\begin{equation}\label{glmn-duality}
\F\cong\sum_{\la\in\Z}L(\hat{gl}_{m|n},\tilde{\La}(\la))\otimes
V^\la_1,
\end{equation}
where $\tilde{\La}(\la)$ is given as follows:
\begin{equation}\label{KWweights}
\tilde{\La}(\la)=
\begin{cases}
\tilde{\epsilon}_1+\cdots+\tilde{\epsilon}_\la+ \tilde{\La}_0,&
\text{for}\ \ 0\le \la\le m,\\
\tilde{\epsilon}_1+\cdots+\tilde{\epsilon}_m+
(\la-m)\tilde{\delta}_1+\tilde{\La}_0,& \text{for}\ \ \la> m,\\
\la\tilde{\delta}_n+ \tilde{\La}_0,& \text{for}\ \ \la<0.
\end{cases}
\end{equation}
We remark that \eqnref{KWweights} is a complete list (up to {\em
essential equivalence} \cite{KW2}) of integrable highest weights
for $\hat{gl}_{m|n}$ of level $1$ with $m\ge 2$ \cite{KW2}. Now we
can write the character of $\F$ as
\begin{equation}\label{fullglmn}
{\rm ch}\F=\prod_{s\in\half+\Z_+}\frac{\prod_{i=1}^m(1+x y_{i}q^s)
(1+x^{-1}y_{i}^{-1}q^{s})}{\prod_{j=1}^n(1-x
z_{j}q^s)(1-x^{-1}z_{j}^{-1}q^s)},
\end{equation}
where again we use the notation $y_i=e^{\tilde{\epsilon}_i}$,
$z_j=e^{\tilde{\delta}_j}$ and $x=e^{\epsilon}$ with $\epsilon\in
(gl_1)^* $ such that $\epsilon(I)=1$. By \propref{combid}
\eqnref{fullglmn} can be written as ($q=e^{\tilde{\delta}}$)
\begin{equation}
\sum_{\la\in\Z}\Bigl{(}\sum_{\mu-\nu=\la}HS_\mu({\bf zq};{\bf
yq})HS_\nu({\bf z^{-1}q};{\bf y^{-1}q})\Bigl{)}x^\la,
\end{equation}
where ${\bf z^{\pm 1}q}=\{z^{\pm 1
}_jq^s|j=1,\cdots,n;s\in\hf\Z_+\}$ and ${\bf y^{\pm 1}q}=\{y^{\pm
1 }_iq^s|i=1,\cdots,m;s\in\hf\Z_+\}$. Thus we have arrived at the
following description of characters by \lemref{independence}.

\begin{thm} Let $\la\in\Z$ and $\tilde{\La}(\la)$ be the $\hat{gl}_{m|n}$-highest weight
defined by \eqnref{KWweights}.  Then
\begin{equation*}
{\rm
ch}L(\hat{gl}_{m|n},\tilde{\La}(\la))=q^{-\frac{|\la|}{2}}\sum_{\mu-\nu=\la}HS_\mu({\bf
zq};{\bf yq})HS_\nu({\bf z^{-1}q};{\bf y^{-1}q}),
\end{equation*}
where $\mu,\nu\in\Z_+$.
\end{thm}

\begin{rem}\label{aux59}
We note that we can derive character formulas in an analogous
fashion for the $\hat{gl}_{\infty}$-modules that appear in the
Fock space decomposition \cite{F} \cite{FKRW} \cite{KR2} (see also
\cite{W1} and \cite{W2} for a rather elegant argument in the
spirits of Howe duality). The resulting formulas will be sums of
products of ordinary Schur functions instead of hook Schur
functions. These formulas agree with the ones obtained in
\cite{KR2} and the one in \cite{AFMO2} in the special case when
$\La=0$. However, when dealing with the $q$-character formulas, we
can also produce a character involving just a sum of Schur
functions.
\end{rem}

\begin{rem}
By \cite{KR1} the representation theory of $\hat{gl}_\infty$ is
closely related to the representation theory of $W_{1+\infty}$,
which is the limit (in an appropriate sense) of the $W$-algebras
$W_N$, as $N\to\infty$ \cite{PRS1} \cite{PRS2}. In particular,
quasi-finite irreducible highest weight representations of latter
can be constructed on a suitable tensor product of quasi-finite
irreducible highest weight representations of a central extension
of $gl_{\infty}\otimes A_n$, where $A_n\cong\C[t]/t^n$. Using
analogous argument it can be shown that the quasi-finite highest
weight irreducible representations of the Lie superalgebra of
differential operators on the super circle, the super
$W_{1+\infty}$ introduced in \cite{MR} (cf. \cite{AFMO} and
\cite{CW3} for definition), can be realized on a suitable tensor
product of quasi-finite highest weight irreducible modules of a
central extension of $gl_{\infty|\infty}\otimes A_n$. Furthermore
each $L(\hgltwo,\La(\la))$ carries a structure of an irreducible
representations of super $W_{1+\infty}$ \cite{CW3}. In particular,
our character formula may be modified to obtain a character
formula for these quasi-finite irreducible super
$W_{1+\infty}$-modules.
\end{rem}

\section{Tensor Product Decomposition}

In this section as another application of \thmref{duality} we will
compute the tensor product decomposition
\begin{equation}\label{tensor}
L(\hat{gl}_{\infty|\infty},\Lambda(\mu))\otimes
L(\hat{gl}_{\infty|\infty},\Lambda(\nu))\cong\sum_{\la}a^{\mu\nu}_\la
L(\hat{gl}_{\infty|\infty},\Lambda(\la)),
\end{equation}
where $\Lambda(\mu)$ and $\Lambda(\nu)$ denote highest weights of
$\hat{gl}_{\infty|\infty}$ of level $l$ and level $r$,
respectively, so that $\mu$ and $\nu$ are generalized partitions
with $l(\mu)\le l$ and $l(\nu)\le r$. The summation $\la$ in
\eqnref{tensor} is over all generalized partitions of length not
exceeding $l+r$ and $\Lambda(\la)$ is viewed as a highest weight
of $\hat{gl}_{\infty|\infty}$ of level $l+r$. We will compute the
coefficients $a^{\mu\nu}_\la$ in terms of the usual
Littlewood-Richardson coefficients (see e.g.~\cite{M}).

To emphasize the dependence of the Fock space $\F$ in
\thmref{duality} on the integer $l$ we will write $\F=\F^l$ and
hence \thmref{duality} reads
\begin{equation*}
\F^l\cong\sum_{\la}L(\hat{gl}_{\infty|\infty},\Lambda(\la))\otimes
V_l^{\la}.
\end{equation*}
Therefore we have
\begin{equation*}
\F^l\otimes\F^r\cong\sum_{\mu,\nu}\Big{(}L(\hat{gl}_{\infty|\infty},\Lambda(\mu))\otimes
L(\hat{gl}_{\infty|\infty},\Lambda(\nu))\Big{)}\otimes
\Big{(}V_l^{\mu}\otimes V_r^\nu\Big{)}.
\end{equation*}
Now $\F^l\otimes\F^r\cong\F^{l+r}$ and hence using
\thmref{duality} again we have
\begin{equation}\label{aux61}
\sum_{\la}L(\hat{gl}_{\infty|\infty},\Lambda(\la))\otimes
V_{l+r}^\la\cong
\sum_{\mu,\nu}\Big{(}L(\hat{gl}_{\infty|\infty},\Lambda(\mu))\otimes
L(\hat{gl}_{\infty|\infty},\Lambda(\nu))\Big{)}\otimes
\Big{(}V_l^{\mu}\otimes V_r^\nu\Big{)}.
\end{equation}
Now suppose that $V_{l+r}^\la$, when regarded as a $gl_l\times
gl_r$-module via the obvious embedding of $gl_l\times gl_r$ into
$gl_{l+r}$, decomposes as
\begin{equation*}
V_{l+r}^\la\cong \sum_{\mu,\nu}b^\la_{\mu\nu}V_l^\mu\otimes
V_r^\nu.
\end{equation*}
This together with \eqnref{tensor} and \eqnref{aux61} give
\begin{equation}\label{aux62}
a^{\mu\nu}_\la=b^\la_{\mu\nu}.
\end{equation}

The duality between the branching coefficients and tensor products
of a general dual pair is well-known \cite{H2}.  We recall that in
\eqnref{aux62} $\mu$, $\nu$ and $\la$ are generalized partitions
subject to constraints on their lengths.

Now \propref{aux31} combined with an analogous argument as the one
given above imply that
\begin{equation}\label{aux63}
\tilde{a}^{\mu\nu}_\la=b^\la_{\mu\nu},
\end{equation}
where here $\mu$, $\nu$, $\la$ are partitions of appropriate
lengths and the $\tilde{a}^{\mu\nu}_\la$'s are the usual
Littlewood-Richardson coefficients.  We remark that there are
combinatorial algorithms to compute these coefficients, the most
well-known probably being the celebrated Littlewood-Richardson
rule (again see e.g.~\cite{M}).

Now for generalized partitions $\mu$, $\nu$ and $\la$ of
appropriate lengths the decomposition $V_{l+r}^\la\cong
\sum_{\mu,\nu}b^\la_{\mu\nu}V_l^\nu\otimes V_r^\mu$ implies that
$V_{l+r}^{\la+d{\bf 1}_{l+r}}\cong
\sum_{\mu,\nu}b^\la_{\mu\nu}V_l^{\mu+d{\bf 1}_l}\otimes
V_r^{\nu+d{\bf 1}_r}$, where here ${\bf 1}_k$ denotes the
$k$-tuple $(1,1,\cdots,1)$ regarded as a partition.  Hence
$b^\la_{\mu\nu}=b^{\la+d{\bf 1}_{l+r}}_{\mu+d{\bf 1}_l,\nu+d{\bf
1}_r}$. Now if we choose a non-negative integer $d$ so that
${\la+d{\bf 1}_{l+r}}$ is a partition, then $b^{\la+d{\bf
1}_{l+r}}_{\mu+d{\bf 1}_l,\nu+d{\bf 1}_r}=\tilde{a}_{\la+d{\bf
1}_{l+r}}^{\mu+d{\bf 1}_l,\nu+d{\bf 1}_r}$ and hence by
\eqnref{aux62} and \eqnref{aux63}
\begin{equation*}
a_\la^{\mu\nu}=\tilde{a}_{\la+d{\bf 1}_{l+r}}^{\mu+d{\bf
1}_l,\nu+d{\bf 1}_r}.
\end{equation*}

From our discussion above we arrive at the following theorem.

\begin{thm} Let $\mu$ and $\nu$ be generalized partitions with
$l(\mu)\le l$ and $l(\nu)\le r$ so that we may regard
$\Lambda(\mu)$ and $\Lambda(\nu)$ as
$\hat{gl}_{\infty|\infty}$-highest weights of level $l$ and $r$,
respectively. Then we have the following decomposition of
$L(\hat{gl}_{\infty|\infty},\Lambda(\mu))\otimes
L(\hat{gl}_{\infty|\infty},\Lambda(\nu))$ into
$\hat{gl}_{\infty|\infty}$-highest weight modules of level $l+r$:
\begin{equation*}
L(\hat{gl}_{\infty|\infty},\Lambda(\mu))\otimes
L(\hat{gl}_{\infty|\infty},\Lambda(\nu))\cong\sum_{(\la,d)}\tilde{a}^{\mu+d{\bf
1}_l,\nu+d{\bf 1}_r}_\la
L(\hat{gl}_{\infty|\infty},\Lambda(\la-d{\bf 1}_{l+r})),
\end{equation*}
where the summation above is over all pairs $(\la,d)$ subject to
the following three conditions:
\begin{itemize}
\item[(i)] $\la$ is a partition of length not exceeding $l+r$
and $d$ a non-negative integer.
\item[(ii)]
$\mu+d{\bf 1}_l$ and $\nu+d{\bf 1}_r$ are partitions.
\item[(iii)]
If $d>0$, then $\la$ is a partition with $\la_{l+r}=0$.
\end{itemize}
Here the coefficients $\tilde{a}^{\mu+d{\bf 1}_l,\nu+d{\bf
1}_r}_\la$ are determined by the tensor product decomposition of
$gl_{k}$-modules $V^{\mu+d{\bf 1}_l}_{k}\otimes V^{{\nu+d{\bf
1}_r}}_{k}\cong\sum_{\la}\tilde{a}^{\mu+d{\bf 1}_l,\nu+d{\bf
1}_r}_\la V^\la_{k}$, where $k\ge l+r$.
\end{thm}

\begin{rem}
Making use of the Howe duality between $\hat{gl}_{\infty}$ and
$gl_l$ mentioned earlier in \remref{aux59}, one derives in a
completely analogous fashion a tensor product decomposition rule
for these $\hat{gl}_{\infty}$-modules that is identical to that
for $\hgltwo$-modules.
\end{rem}

\end{document}